\documentclass[12pt]{article}
\usepackage{amssymb,amsmath,latexsym}

{\bf}{\it}
\newtheorem{definition}{Definition}{\bf}{\it}
{\bf}{\rm}
\newtheorem{lemma}{Lemma}{\bf}{\it}
\newtheorem{theorem}{Theorem}{\bf}{\it}
{\bf}{\it}
{\bf}{\it}
{\bf}{\rm}
\def\al{\alpha}

\def\ap{\rightarrow}
\def\A{{\cal A}}

\def\eb{{\bf e}}

\def\f{{\bf f}}

\def\F{{\cal F}}

\def\Pi{P^{\infty}}

\renewcommand{\Re}{\mathbb{R}}
\newcommand{\R}{\mathbb{R}}
\newcommand{\Rn}{\R^n}
\newcommand{\vol}{\mathrm{Vol}}
\newcommand{\bO}{\mathbf{O}}

\def\R{{\bf R}}

\def\t{\tau}
\def\th{\theta}

\def\u{{\bf u}}

\def\v{{\bf v}}
\def\vcd{\mathrm{VCDIM}}
\def\w{{\bf w}}
\def\x{{\bf x}}

\def\xm1{x_{m+1}}
\def\xim1{\xi_{m+1}}
\def\y{{\bf y}}

\def\seq{\subseteq}
\newcommand{\qed}{\hfill $\Box$}
\newcommand{\dia}{\hfill $\diamond$}
\newcommand{\rd}{{\mathrm{\underline{rd}}}}

\newcommand{\bc}{\begin{center}}
\newcommand{\ec}{\end{center}}
\newcommand{\be}{\begin{equation}}
\newcommand{\ee}{\end{equation}}
\newcommand{\bd}{\begin{displaymath}}
\newcommand{\ed}{\end{displaymath}}
\newcommand{\ba}{\begin{array}}
\newcommand{\ea}{\end{array}}
\newcommand{\ben}{\begin{enumerate}}
\newcommand{\een}{\end{enumerate}}
\newcommand{\bit}{\begin{itemize}}
\newcommand{\eit}{\end{itemize}}
\newcommand{\beq}{\begin{eqnarray}}
\newcommand{\eeq}{\end{eqnarray}}

\renewcommand{\Xi}{\mbox{$X^{\infty}$}}

\renewcommand{\and}{\mbox{$\wedge$}}

\def\bi{\{0,1\}}

\begin{document}
\thispagestyle{empty}

\title{An Improved Bound on the VC-Dimension of\\
Neural Networks with Polynomial Activation Functions}

\author{
\begin{tabular}{cc}
M.\ Vidyasagar & J.\ Maurice Rojas \\
Advanced Technology Centre & Department of Mathematics \\
 Tata Consultancy Services & Texas A\&M University \\
1-2-10 S.\ P.\ Road & College Station, TX \ 77843-3368 \\
INDIA 500 003 & U.S.A.\\
{\tt sagar@atc.tcs.co.in} & {\tt rojas@math.tamu.edu} 
\end{tabular}
}

\date{\today}

\maketitle
\begin{abstract}
We derive an improved upper bound for the VC-dimension of
neural networks with polynomial activation functions.
This improved bound is based on a result of Rojas \cite{Rojas00}
on the number of connected components of a semi-algebraic set.
\end{abstract}

\newpage

\baselineskip 18pt

\section{Introduction}\label{sec:1}
We examine neural networks with polynomial activation functions.
The specific architecture of the neural networks is described in detail in
the next section.
Such neural networks have been the subject of active investigation for several
years, since powerful tools from algebraic geometry can be brought into
play in analyzing the VC-dimension of such networks.
Perhaps \cite{GJ93} was the first paper to connect these two subjects.
For several years (see for example \cite{KM97}) it has been known that
every bound on the number of connected components of a semi-algebraic set  
can be readily translated into a corresponding bound on the VC-dimension
of a neural network architecture.
Practically all of the known bounds on the VC-dimension of neural networks
with polynomial activation bounds make use of a classical result
discovered by Milnor, Oleinik, Petrovsky, and Thom 
\cite{op,milnor,thom}.\footnote{Actually, this result bounds 
the sum of the Betti numbers of a semi-algebraic set, and this  
quantity is always at least as large as the number of connected components.
In practice, one usually only needs an upper bound on the number of connected 
components.} This bound, while easy to use, is usually much larger than 
necessary, since it only uses coarse information about the underlying 
set such as the number of variables and the maximum 
degree of the input polynomials. 
More recently, sharper bounds using more refined data from the input  
polynomials have been discovered. In the present note, we use a result due to 
Rojas \cite{Rojas00} that is particularly well-suited to neural networks with 
polynomial activation functions.
The present bound is, {\em in all cases}, sharper than the
earlier bound of Goldberg and Jerrum \cite{GJ93}.
Moreover, it is intuitively appealing, as the improvement can be quantified
as the relative entropy of two probability vectors, whose dimension equals
the number of layers in the neural network.
This shows that the problem of bounding the VC-dimension of a neural network
architecture continues to be interesting, and that we should strive
to derive even tighter upper bounds.

Our main result is stated in theorem \ref{thm:new} of section \ref{sec:4}. 
The recent semi-algebraic bound it is based on is stated as theorem 
\ref{thm:Rojas} of section \ref{sec:3}. However, let us first review a bit of 
background and some of the earlier bounds. 

\section{Known Results}\label{sec:2}

The following definition of the VC-dimension is standard; see
for example the books by Vapnik \cite{Vapnik95} or Vidyasagar \cite{MV97}.
\begin{definition}\label{def:VC-dim}
Suppose $X$ is a set and $\F$ is a collection of $\{0,1\}$-valued
functions on $X$.
A set $S = \{ x_1 , \ldots , x_n \} \seq X$ is said to be {\bf shattered}
by $\F$ if each of the $2^n$ functions mapping $S$ into $\bi$ is the
restriction to $S$ of some function in $\F$.
The {\bf Vapnik-Chervonenkis (VC)-dimension} of $\F$ is the largest
integer $n$ such that there exists a set of cardinality $n$ that is
shattered by $\F$. \dia  
\end{definition}
By identifying a $\{0,1\}$-valued function with its support set,
it is also possible to speak of the VC-dimension of a collection of
sets.
In the sequel, we shall use both notions interchangeably.

Following by now familiar approaches, we view a {\em neural network\/}
as a verifier of formulas.
Specifically, let $\w \in \Re^k$ denote the ``weight vector'' or the
vector adjustable parameters in a neural networks.
A neural network with input space $X \seq \Re^N$ and  
weight vector $\w$ evaluates a logical proposition $\phi(\x,\w)$ which is
a Boolean combination of $s$ atomic expressions of the form 
$\t_i(\x,\w)\!=\!0$ or $\t_j(\x,\w)\!>\!0$. Letting $1$ 
(resp.\ $0$) denote ``true'' (resp.\ ``false''), we can thus think 
of $\phi$ as a function from $\Re^{k+N}$ to $\{0,1\}$. 
So for each weight vector $\w$, define
\bd
A_{\w} := \{ \x \in \Re^N : \phi(\w,\x) = 1 \} .
\ed
The objective is to obtain an upper bound on the VC-dimension of
the collection of sets $\A := \{ A_{\w} \; | \; \w \in \Re^k\}$ or, 
equivalently, the VC-dimension of the collection of $\{0,1\}$-valued functions
$\Phi := \{ \phi(\w,\cdot) \; | \;  \w \in \Re^k \}$.

To state the result, we need one final bit of notation: 
Let $\x_1 , \ldots , \x_v
\in \Re^N$, and suppose $sv \geq k$.
{}From the $sv$ polynomials $\t_j(\cdot,\x_i)$ 
determined by all $(i,j)\!\in\!\{1,\ldots,v\}\times \{1,\ldots,s\}$, 
choose $r \leq k$ polynomials, and label them 
$\th_1(\cdot), \ldots , \th_r(\cdot) : \Re^k \ap \Re$.
Define
\bd
\f(\w) := [\th_1(\w) \ldots \th_r(\w)] \in \Re^r .
\ed
Finally, let $B$ denote the maximum number of connected 
components of any pre-image $f^{-1}(y)$ with $y\!\in\!\R^r$, 
for any choice of $r$ and $\theta_1,\ldots,\theta_r$ as above. 
With the above set-up, the following result is proved in \cite{KM97}.
For further background on our setting below see \cite{KM97} or \cite{MV97}, 
p.\ 329.

\begin{theorem}\label{thm:KM}
Following the notation above, assume further that 
restrict to those $y\!\in\!\R^r$ that are regular values of 
$f$. Then  
\[\label{eq:21}
\vcd(\Phi) \leq 2 \lg B + 2k \; \lg (2es). 
\]
\end{theorem}
That $B$ is in fact finite and admits an explicit upper bound 
is obtained by appealing to the aforementioned classical result of 
Milnor, Oleinik, Petrovsky, and Thom \cite{op,milnor,thom}, which we 
now state as follows:
\begin{lemma}\label{lemma:Milnor}
Suppose $\th_1 , \ldots , \th_r$ are polynomials in $k$ variables, with
degree no larger than $d$.
Then whenever $\y$ is a regular value of $\f$ as defined above,
the preimage $\f^{-1}(\y)$ contains no more than $d(2d)^{k-1}$
connected components.
\end{lemma}
Note that Milnor actually proves the theorem in the case where $\y = {\bf 0}$,
but we can clearly perturb the constant terms of the $\theta_i$ 
to enforce this assumption. 
If we replace the quantity $d(2d)^{k-1}$ by the larger number $(2d)^k$
and substitute $B = (2d)^k$ into the upper bound (\ref{eq:21}),
we get the following result.
\begin{theorem}\label{thm:GJ}
Suppose $\phi$ is a Boolean formula involving a total of $s$ polynomial 
equalities and inequalities, where each polynomial has degree no larger than 
$d$ with respect to $\w$. Then
$\vcd(\Phi) \leq 2k\; \lg (4eds)$. 
\end{theorem}

The above result is the same as that derived in \cite{GJ93}.
It should be noted, however, that Goldberg and Jerrum actually consider
neural networks with {\em piecewise polynomial\/} activation functions.
With more elaborate notation, their results can be derived as special
cases of Theorem \ref{thm:KM}.

Theorem \ref{thm:KM} shows the importance of deriving {\em tight\/}
upper bounds on the number of connected components of a semi-algebraic set.
This is a long-standing problem in real algebraic geometry
that has received considerable attention from the research community.
It is obvious from the bound (\ref{eq:21}) that any improvement over
Milnor's upper bound translates {\em directly\/} into a corresponding
improvement in the estimate of the VC-dimension of a neural network
architecture with polynomial activation functions.
This leads us to the next topic.

\section{Improved Upper Bound on the Number of Connected
Components}\label{sec:3}

\setcounter{equation}{0}

In \cite{Rojas00}, an improvement is provided over Milnor's bound.
To state this improved result, a bit of notation is introduced.

Let $\Delta_n$ denote the standard $n$-simplex in $\Rn$, with 
vertices the standard basis vectors {\em and} the origin. 
Note that 
\bd
d\Delta_n := \left\{ (x_1 , \ldots , x_n)\!\in\!\Re^n 
\; | \;  x_i \geq 0 \text{ for all } i  \text{ and } 
\sum_{i=1}^n x_i \leq d \right\} .
\ed
Let $\vol_n(\cdot)$ denote the renormalization of the usual volume in 
$\Rn$ satisfying $\vol_n(\Delta_n) = 1$. (Since the usual $n$-dimensional 
volume is multiplicative for orthogonal subspaces, it is easy to 
prove by induction that $\vol_n$ is just $n!$ times the usual 
$n$-dimensional volume.) 
\begin{theorem}\label{thm:Rojas}
Suppose $\t_1, \ldots , \t_r$ are polynomials in $(w_1,\ldots,w_k)$, 
and let $\eb_1 , \ldots , \eb_k$ and $\bO$ denote the standard basis vectors
and the origin of $\Re^k$. Also, let $V$ denote the convex hull of the 
union of $\{\bO,e_1,\ldots,e_k\}$ with the set of all 
$(i_1 , \ldots , i_k)$ such that
$w_1^{i_1} \ldots w_k^{i_k}$ is a monomial of some $\th_j(\cdot)$.  Then
\[\label{eq:31}
B \leq 2^k \vol_k(V). 
\]
\end{theorem}
In the special case where {\em every\/} $k$-tuple with $\sum_{j=1}^k
i_j \leq d$ occurs in $V$, we recover the (adjusted) Milnor bound
$(2d)^k$. However, the whole point of the preceding refined 
bound is that there are many instances where the input 
polynomial are far more sparse, and this can be exploited. 

\section{Improved Upper Bound on the VC-Dimension}\label{sec:4}

\setcounter{equation}{0}

In this section, we derive an improved upper bound on the VC-dimension
of neural networks with polynomial activation functions.
The improved bound is a direct consequence of coupling Theorems
\ref{thm:KM} and \ref{thm:Rojas}.

Let us begin by describing the class of neural networks under study.
It is assumed that the network has $N$ real inputs denoted by $x_1 ,
\ldots , x_N$. 
There are $l$ levels in the network, and at level $i$ there are 
$q_i$ output neurons; however, at the output layer (level $l$)
there is only a single neuron (see below).
Let $k_i$ denote the number of adjustable parameters, or ``weights,''
at level $i$, and let $k = \sum_{i=1}^l k_i$ denote the total
number of adjustable parameters.
Let $\w_i := (w_{i,1}, \ldots , w_{i,k_i})$ denote the weight vector at
level $i$, and $\w = (\w_1 \ldots \w_l)$ denote the total weight vector.
The input-output relationship of each neuron at level $i$ is of the form
\bd
y_{i,j} = \t_{i,j}(\w_i, y_{i-1,1} , \ldots , y_{i-1,q_{i-1}} ) , \;
j = 1, \ldots , q_i .
\ed
where $y_{i,j}$ is the output of neuron $j$ at level $i$, and
$\t_{i,j}$ is a polynomial of degree no larger than $\al_i$ in the
components of the weight vector $\w_i$, and no larger $\beta_i$
in the components of the vectors $y_{i-1,j}$.
At the final layer, there is a simple perceptron device following the
polynomial activation function.

With this class of neural networks, it is clear that the output will equal
one if and only if a polynomial inequality of the form
\bd
y_l(\w,\x) \geq 0 ,
\ed
is satisfied, where $\w$ is the weight vector and $\x = (x_1 \ldots x_N)$
is the input vector.
Thus we can apply Theorem \ref{thm:KM} with $s = 1$.
The issue now is to determine the number of connected components $B$
of the semi-algebraic set defined by $y_l(\w,\x) = \y$.

Now we are in a position to state the main result.
To facilitate the statement, we introduce a bit of notation.
Define
\bd
d_l = \al_l , \; d_{l-1} = \al_{l-1} \beta_l, \ldots , d_i = \al_i
\prod_{j=i+1}^l \beta_j ,
\; i = 1, \ldots , l-1 .
\ed
Recall that $k_i$ denotes the number of adjustable parameters at level $i$,
and that $k$ denotes the total number of adjustable parameters.
Define the probability vectors
\bd
\v := (k_1/k \ldots k_l/k) , \; \u := (d_1/d \ldots , d_l/d), 
\ed
and define the ``binary'' relative entropy $H(\v|\u)$ as
\bd
H(\v|\u) := \sum^l_{i=1}v_i\lg\left(\frac{v_i}{u_i}\right)\!=\!
\frac{1}{k}\sum_{i=1}^l  k_i \lg\left(\frac{dk_i}{kd_i}\right).
\ed
Note that the above is the same as the conventional relative entropy
of two probability vectors, except that we use base-$2$ logarithms instead
of natural logarithms.
Following standard convention, we take $0 \lg (0/0) = 0$.
\begin{theorem}\label{thm:new}
With the above notation, we have
\beq
B & \leq & 2^k k! \prod_{i=1}^l \frac{ d_i^{k_i} }{ k_i !} \label{eq:41}. 
\\ & \leq & \left(\frac{2d}{e^{7/8}}\right)^k 2^{-kH(\v|\u)} . 
\label{eq:42}
\eeq
where $d\!:=\!\sum^l_{i=1}d_i$ and we assume $k_1,\ldots,k_l\!\geq\!2$ 
in the last inequality. 
Consequently, when $k_1,\ldots,k_n\!\geq\!2$, the VC-dimension of the neural 
network architecture is bounded above by
\[\label{eq:43} 2k(\lg (4ed) - H(\v|\u)).  \]
\end{theorem}

\noindent
{\bf Remark} The above theorem shows that the reduction in the 
VC-dimension estimate over that of Theorem \ref{thm:GJ} is precisely $2k$ 
times the (binary) relative entropy of the two probability vectors $v$ 
and $u$ defined above. Thus if $k_i/k = d_i/d$ for all $i$, there will not 
be any reduction at all.  In general, the fraction by which the 
older VC-dimension estimate is reduced is precisely the 
ratio $H(\v|\u)/(\lg(4ed))$. Note also that the assumption that 
there are at least $2$ adjustable parameters at each levels is 
a reasonably mild assumption. $\diamond$  

\noindent
{\bf Proof of Theorem \ref{thm:new}:}
The proof depends on a careful book-keeping of the degree of $y_l(\w,\x)$
with respect to the various components of $\w$.
{}From the architecture of the neural network, it is clear that at the first
level, each of the $y_{1,j}$ is a polynomial in the components of $\w_1$
of degree no larger than $\al_1$.
At the second level, each of the $y_{2,j}$ is a polynomial, whose monomials
are of (combined) degree no larger than $\al_2$ in the components of $\w_2$,
and of (combined) degree no larger than $\beta_2 \al_1$ in the components
of $\w_1$.
Thus, while each $y_{2,j}$ could have a total degree of $\al_2 + \beta_2 \al_1$
in the components of $\w_1$ and $\w_2$, the total degree of the monomial
terms involving the components of $\w_1$ does not exceed $\beta_2 \al_1$, while
the total degree of the monomial terms involving the components of $\w_2$
does not exceed $\al_2$.
A simple argument by induction then tells us that at the output layer (level 
$l$), the single output $y_l$ is a polynomial whose monomials have total 
degree no larger than $d_l = \al_l$ in the components of $\w_l$,
no larger than $d_{l-1} = \beta_l \al_{l-1}$ in the components of $\w_{l-1}$,
and so on.
With the $d_i$'s defined as above, the components of each $\w_i$ appear
with total degree no larger than $d_i$.
Thus the total degree of $y$ could be as large as $d$,
but the monomial terms involving the components of $\w_i$ have total degree
no larger than $d_i$.
So the set $V$ defined in Theorem \ref{thm:Rojas} satisfies the following
containment:
\bd
V \seq \prod_{i=1}^l S_{d_i}^{k_i} .
\ed
Because of this containment, it follows that 
\bd
\vol_n(V) \leq k!\prod_{i=1}^l \frac{ d_i^{k_i} }{k_i !}. 
\ed 
Combining this with the bound (\ref{eq:31}) establishes the first estimate
(\ref{eq:41}).

To prove the second estimate, we use Stirling's approximation. 
In particular, \cite[ex.\ 20, pg.\ 200]{rudin} tells us that 
for all $t\!\in\!\{2,3,4,\ldots\}$, we have 
\[ e^{7/8} (t/e)^t\sqrt{t} < t! < e (t/e)^t \sqrt{t}. \] 
Consequently, we easily obtain 
\[ \frac{k!}{k_1!\cdots k_l!} < e^{1-\frac{7}{8}k} \frac{k^k}
{k^{k_1}_1\cdots k^{k_l}_l\sqrt{k_1\cdots k_l}}.\] 
Dropping the square root term on the bottom can of course  
be done, and then an elementary calculation yields 
$2^k k! \prod_{i=1}^l \frac{ d_i^{k_i} }{ k_i !}\!\leq\!
\left(\frac{2d}{e^{7/8}}\right)^k 2^{-kH(\v|\u)}$, provided 
$k_1,\ldots,k_n\!\geq\!2$. 

The VC-dimension estimate (\ref{eq:43}) now follows readily from
Theorem \ref{thm:KM}.  \qed

\section{Numerical Example}\label{sec:5}

Consider a network with four inputs, five hidden-layer
neurons at the first level and an output-layer neuron.
As is common, let us suppose that $\al_i = 1$ for all $i$.
This means that all the adjustable parameters {\em enter linearly\/}
into the corresponding activation function.
Suppose $\beta_1 = 2, \beta_2 = 3$.
This means that the hidden-layer neurons have quadratic activation functions,
whereas the output-layer neuron has a cubic activation function.
It remains to specify the integers $k_1$ and $k_2$, representing
the number of adjustable parameters.
Let us assume that practically all of the monomial terms are present in
each neural characteristic.
Thus it is reasonable to assume $k_1 = 50, k_2 = 20$.
Finally, $d_1 = 3, d+2 = 1$.
With these figures, one has
\bd
\v = (5/7 , 2/7) \ , \  \u = (0.25 , 0.75) ,
\ed
\bd
H(\v|\u) \approx 0.684033, \; \lg(4ed) \approx 5.4427, \;
\frac{H(\v|\u)}{\lg(4ed)} \approx 0.12567 .
\ed
Thus, in this case, the improved bound is roughly 12.5\% sharper. 

\section{Conclusions}\label{sec:6}

\end{document}